\documentclass{article}

\usepackage{times}
\usepackage{latexsym}
\usepackage{amsmath}
\usepackage{amssymb}
\usepackage{graphicx}

\newcommand{\nc}{\newcommand}
\nc{\rnc}{\renewcommand}
\nc{\nev}{\newenvironment}

\makeatletter
\rnc{\@seccntformat}[1]{{\normalfont\bfseries{\csname the#1\endcsname}\hspace{0.5em}}}
\rnc{\section}{\@startsection
        {section}%
        {1}%
        {0mm}%
        {-\baselineskip}%
        {0.5\baselineskip}%
        {\normalfont\normalsize\bfseries\centering}%
}
\rnc{\subsection}{\@startsection
        {subsection}%
        {2}%
        {0mm}%
        {-0.5\baselineskip}%
        {-0.5em}%
        {\normalfont\normalsize\bfseries}%
}
\makeatother

\nc{\sus}[1]{\subsection*{#1.}}

\newcounter{theo}[section]
\rnc{\thetheo}{\arabic{section}.\arabic{theo}}
 
\nev{defform}[2]{\refstepcounter{theo}\begin{list}{}{%
      \setlength{\labelwidth}{0em}%
      \setlength{\leftmargin}{0em}%
      \setlength{\itemindent}{\labelsep}%
      \setlength{\listparindent}{1.5em}
      \setlength{\parsep}{0em}}%
      \item[\textbf{#1 \thetheo#2.}]}%
    {\end{list}}

\nev{satzform}[2]{\begin{defform}{#1}{#2}\itshape}{\end{defform}}

\nev{ndef}[2]{\begin{defform}{#1}{ (#2)}}{\end{defform}}
\nev{nsatz}[2]{\begin{defform}{#1}{ (#2)}\itshape}{\end{defform}}

\nev{theo}{\begin{satzform}{Theorem}{}}{\end{satzform}}
\nev{lem}{\begin{satzform}{Lemma}{}}{\end{satzform}}
\nev{cor}{\begin{satzform}{Corollary}{}}{\end{satzform}}
\nev{prop}{\begin{satzform}{Proposition}{}}{\end{satzform}}
\nev{fact}{\begin{satzform}{Fact}{}}{\end{satzform}}
\nev{ex}{\begin{defform}{Example}{}}{\end{defform}}
\nev{defn}{\begin{defform}{Definition}{}}{\end{defform}}
\nev{rem}{\begin{defform}{Remark}{}}{\end{defform}}

\nc{\proof}{\medskip\noindent\textsc{Proof: }}
\nc{\proofend}{\unskip\nobreak\hfill\raisebox{-.067em}{\large$\Box$}\vspace{\topsep}\par}

\rnc{\labelenumi}{(\arabic{enumi})}
\rnc{\labelitemi}{\text{--}}
\rnc{\phi}{\varphi}

\nc{\bigmid}{\;\big\|\;}
\nc{\Bigmid}{\;\Big\|\;}

\nc{\val}{\textnormal{val}}
\rnc{\int}{\textnormal{int}}
\rnc{\max}{\textnormal{max}}
\rnc{\min}{\textnormal{min}}
\nc{\ad}{\textnormal{ad}}
\nc{\tw}{\textnormal{tw}}
\nc{\ltw}{\textnormal{ltw}}
\nc{\opt}{\textnormal{opt}}
\nc{\var}{\textnormal{var}}
\nc{\Pow}{\textnormal{Pow}}
\nc{\comp}{\textnormal{comp}}
\nc{\valu}{\textnormal{value}}
\nc{\weight}{\textnormal{weight}}
\nc{\dom}{\textnormal{dom}}

\nc{\true}{\textsc{True}}
\nc{\false}{\textsc{False}}

\nc{\maxh}{{\normalfont\scshape Maximum $H$-Matching}}
\nc{\maxind}{{\normalfont\scshape Maximum Independent Set}}
\nc{\minvc}{{\normalfont\scshape Minimum Vertex Cover}}
\nc{\minds}{{\normalfont\scshape Minimum Dominating Set}}

\begin{document}
\title{\large\bfseries Local tree-width, excluded minors, and approximation algorithms}
\author{\normalsize Martin Grohe\\
\small Institut f\"ur Mathematische Logik, Eckerstr.~1, 79104
Freiburg, Germany\\
\small Email: grohe@logik.mathematik.uni-freiburg.de}
\date{\normalsize January 2000}
\maketitle

\begin{abstract}
  The \emph{local tree-width} of a graph $G=(V,E)$ is the function
  $\ltw^G:\mathbb N\to\mathbb N$ that associates with every
  $r\in\mathbb N$ the maximal tree-width of an $r$-neighborhood in
  $G$. Our main graph theoretic result is a decomposition theorem for
  graphs with excluded minors that essentially says that such graphs
  can be decomposed into trees of graphs of bounded local tree-width.
    
  As an application of this theorem, we show that a number of
  combinatorial optimization problems, such as \minvc, \minds, and
  \maxind\ have a polynomial time approximation scheme when restricted
  to a class of graphs with an excluded minor.
\end{abstract}

\section{Introduction}
\emph{Tree-width}, measuring the similarity of a graph with a tree,
has turned out a to be an important notion both in structural graph
theory and in the theory of graph algorithms. It is well known that
planar graphs may have arbitrarily large tree-width. However, for
every fixed $d$ the class of planar graphs of diameter at most $d$ has
bounded tree-width. In other words, the tree-width of a planar graph
can be bounded by a function of the diameter of the graph. This makes
it possible to decompose planar graps into families of graphs of
small tree-width in an orderly way. Such decompositions of planar
graphs, better known under the name \emph{outerplanar decompositions},
have been explored in various algorithmic settings
\cite{bak94,epp95,khamot96,frigro99}.  The main ideas go back to a
fundamental article of Baker \cite{bak94} on approximation algorithms
on planar graphs.

The \emph{local tree-width} of a graph $G=(V,E)$ is the
function $\ltw^G:\mathbb N\to\mathbb N$ that associates with every
$r\in\mathbb N$ the maximal tree-width of an $r$-neighborhood in $G$. More formally,
we define the \emph{$r$-neighborhood} $N_r(v)$ of a vertex $v\in V$ to be the set of
all $w\in V$ of distance at most $r$ from $v$, and we let $\langle
N_r(v)\rangle$ denote the subgraph induced by $G$ on $N_r(v)$. Then, denoting
the tree-width of a graph $H$ by $\tw(H)$, we let 
\[
\ltw^G(r):=\max\Big\{\tw\big(\langle N_r(v)\rangle\big)\;\Big|\;v\in
V\Big\}.
\]
We are mainly interested in classes of graphs of \emph{bounded local
  tree-width}, that is, classes $\mathcal C$ for which there is a
function $f:\mathbb N\to\mathbb N$ such that for all $G\in\mathcal C$
and $r\in\mathbb N$ we have $\ltw^G(r)\le f(r)$. The class of
planar graphs is an example. It has been observed by Eppstein
\cite{epp95} that if a class $\mathcal C$ is closed under taking
  minors and has bounded local tree-width
(Eppstein calls this the ``diameter-treewidth property''), then the
graphs in $\mathcal C$ admit a decomposition into graphs of small
tree-width in the style of the outerplanar decomposition of planar
graphs, and the planar-graph algorithms based on this decomposition
generalize to graphs in $\mathcal C$.  Eppstein
gave a nice characterization of such classes; he proved that a minor
closed class $\mathcal C$ of graphs has bounded local tree-width if,
and only if, it does not contain all apex graphs.

The main graph-theoretic result of this paper, Theorem
\ref{theo:minorclosed}, can be phrased as follows: Let $\mathcal C$ be
a minor closed class of graphs that does not contain all graphs. Then
all graphs in $\mathcal C$ can be decomposed into a tree of graphs
that, after removing a bounded number of vertices, have bounded local
tree-width. (Of course the converse is also true, but trivial: If
$\mathcal C$ is a minor closed class of graphs such that every graph
in $\mathcal C$ admits such a decomposition, then $\mathcal C$ is not
the class of all graphs.) The proof of this result is based on a deep
structural characterization of graphs with excluded minors due to
Robertson and Seymour \cite{GMXVI}.

We defer the precise technical statement of our decomposition theorem
to Section~\ref{sec:minclosed} and turn to its applications now. In
this paper, we focus on approximation algorithms. But let me
mention that the theorem can also be used to re-prove a
result of Alon, Seymour, and Thomas \cite{aloseytho90a} that graphs $G$
with an excluded minor have tree-width $O(\sqrt{|G|})$ (see
Section~\ref{sec:ast}).\footnote{We have observed this in discussions
  with Reinhard Diestel and Daniela K\"uhn.} 

Actually, the main result of Alon, Seymour, and Thomas's article is a
separator theorem for graphs with an excluded minor, generalizing a
well-known separator theorem due to Lipton and Tarjan \cite{liptar79}
for planar graphs. These separator theorems have numerous algorithmic
applications, among them a  polynomial time approximation scheme (PTAS)
for the \maxind\ problem on planar graphs \cite{liptar80} and, more
generally, classes of graphs with an excluded minor \cite{aloseytho90b}.

A different approach to approximation algorithms on planar graphs is
Baker's \cite{bak94} technique based on the outerplanar decomposition.
It does not only give another PTAS for \maxind, but also for other
problems, such as \minds, to which the technique based on the
separator theorem does not apply. 

We can use our decomposition theorem to extend Baker's approach to
arbitrary classes of graphs with an excluded minor. Our purpose here
is to explain the technique and not to give an extensive list of
problems to which it applies. We show in detail how to get a PTAS for
\minvc\ on classes of graphs with an excluded minor and then explain
how this PTAS has to be modified to solve the problems \minds\ and
\maxind. It should be no problem for the reader to apply the same
technique to other optimization problems.

The paper is organized as follows: In Section~\ref{sec:prel} we fix
our terminology and recall a few basic facts about tree-decompositions
of graphs. Local tree-width is introduced in Section~\ref{sec:ltw}. In
Section~\ref{sec:minclosed}, we prove our decomposition theorem for
classes of graphs with an excluded minor. Approximation algorithms
are discussed in Section~\ref{sec:aa}, and in Section~\ref{sec:appl}
we briefly explain two other applications of the decomposition theorem.

\section{Preliminaries}\label{sec:prel}
The vertex set of a graph $G$ is denoted by $V^G$, the edge set by
$E^G$. Graphs are always assumed to be finite, simple, and undirected. We write
$vw\in E^G$ to denote that there is an edge from $v$ to $w$. For a
subset $X\subseteq V^G$, we let $\langle X\rangle^G$ denote the
induced subgraph of $G$ with vertex set $X$. We let $G\setminus
X:=\langle V^G\setminus X\rangle^G$. For graphs $G$ and $H$, we let
$G\cup H:=(V^G\cup V^H,E^G\cup E^H)$.  We often omit superscripts $^G$
if $G$ is clear from the context.

$K_n$ denotes the complete graph with $n$ vertices, and for an
arbitrary set $X$, $K_X$ denotes the complete graph with vertex set
$X$. A vertex set $X\subseteq V^G$ in a graph $G$ is a \emph{clique}
if $K_X\subseteq G$. The \emph{clique number} $\omega(G)$ of a graph
$G$ is the maximal size of a clique in $G$. For a class $\mathcal C$ of
graphs, we let $\omega(\mathcal C)$ be the maximum of the clique
numbers of all graphs in $\mathcal C$, or $\infty$, if this maximum
does not exist.

Note that if $\mathcal C$ is closed under taking subgraphs and is not the
class of all graphs, then $\omega(\mathcal C)$ is finite.

\sus{Graph minors}
A \emph{minor} of a graph $G$ is a graph $H$ that can be obtained from
a subgraph of $G$ by contracting edges; we write $H\preceq G$ to
denote that $H$ is a minor of $G$. 

Note that $H\preceq G$ if, and only if, there is a mapping
$h:V^H\rightarrow\Pow(V^G)$ such that $\langle h(x)\rangle^G$ is a
connected subgraph of $G$ for all $x\in V^H$, $h(x)\cap
h(y)=\emptyset$ for $x\neq y\in V^H$, and for every edge $xy\in
E^H$ there exists an edge $uv\in E^G$ such that $u\in h(x),v\in
h(y)$. We say that the mapping $h$ \emph{witnesses} $H\preceq G$ and
write $h:H\preceq G$.

A class $\mathcal C$ is \emph{minor closed} if, and only if, for all
$G\in\mathcal C$ and $H\preceq G$ we have $H\in\mathcal C$. We call
$\mathcal C$ non-trivial if it is not the class of all graphs.

A class $\mathcal C$ is \emph{$H$-free} if $H\not\preceq G$ for all
$G\in\mathcal C$. We then call $H$ an \emph{excluded minor} for $\mathcal
C$. Note that a class $\mathcal C$ of graphs has an
excluded minor if, and only if, there is an $n\ge 1$ such that
$\mathcal C$ is $K_n$-free. Furthermore, this is equivalent to saying
that $\mathcal C$ is contained in some non-trivial minor closed class
of graphs.

Robertson and Seymour's \cite{GMXX} \emph{Graph Minor Theorem} states
that for every  minor closed class $\mathcal C$ of graphs there is a
  finite set $\mathcal F$ of graphs such that
\[
\mathcal C=\{G\mid\forall H\in\mathcal F:\;H\not\preceq G\}.
\]
For a nice introduction to graph minor theory we refer
the reader to the last chapter of \cite{die00}, a recent survey is
\cite{rtho99}.

\sus{Tree-decompositions}\label{ss:treedec}
In this paper, we assume trees to be directed from the root to the
leaves. If $tu\in E^T$ we call $u$ a \emph{child} of $t$ and $t$
the \emph{parent} of $u$. The root of a tree $T$ is always denoted by $r^T$. 

A \emph{tree-decomposition} of a graph $G$ is a pair $(T,(B_t)_{t\in
  V^T})$, where $T$ is a tree and $(B_t)_{t\in V^T}$ a family of subsets
of $V^G$ such that $\bigcup_{t\in V^T}\langle B_t\rangle^G=G$ and for
every $v\in V^G$ the set $\{t\mid v\in B_t\}$ is connected.  
The sets $B_t$ are
called the \emph{blocks} of the decomposition.  The \emph{width} of
$(T,(B_t)_{t\in V^T})$ is the number $\max\{\|{B_t}\|\mid t\in V^T\}-1$.
The \emph{tree-width} of $G$, denoted by $\tw(G)$, is the minimal
width of a tree-decomposition of $G$.

The following lemma collects a few simple and well-known facts about
tree-de\-compositions:

\begin{lem}\label{lem:treedec}
\begin{enumerate}
\item
Let $(T,(B_t)_{t\in V^T})$ be a tree-decomposition of a graph $G$ and
$X\subseteq V^G$ a clique. Then there is a $t\in V^T$ such that
$X\subseteq B_t$.

\item
Let $G,H$ be graphs such that
$V^G\cap V^H$ is a clique in both $G$ and $H$. Then $\tw(G\cup
H)=\max\{\tw(G),\tw(H)\}$.

\item
Let $G$ be a graph and $X\subseteq V^G$. Then $\tw(G)\le\tw(G\setminus
X)+|X|$.

\item
Let $G,H$ be graphs such that $H\preceq G$. Then $\tw(H)\le\tw(G)$.
\end{enumerate}
\end{lem}

Throughout this paper, for a tree-decomposition $(T,(B_t)_{t\in V^T})$
and $t\in T\setminus\{r^T\}$ with parent $s$ we let $A_t:=B_t\cap
B_s$. We let $A_{r^T}:=\emptyset$.
 
The \emph{adhesion} of
  $(T,(B_t)_{t\in V^T})$ is the number
\[
\ad(T,(B_t)_{t\in V^T}):=\max\{\|A_t\|\mid t\in V^T\}.
\]
The \emph{torso} of $(T,(B_t)_{t\in V^T})$ at $t\in V^T$ is the subgraph
\[
[B_t]:=\langle B_t\rangle^G\cup K_{A_t}\cup\bigcup_{u\text{ child of }t}K_{A_u},
\]
or equivalently, the subgraph with vertex set $B_t$ in which two vertices are
adjacent if, and only if, either they are adjacent in $G$ or they both
belong to a block $B_u$, where $u\neq t$.
$(T,(B_t)_{t\in V^T})$ is a tree-decomposition of $G$ \emph{over} a class
$\mathcal B$ of graphs if all its torsos belong to $\mathcal B$. 

Note that the adhesion of a tree-decomposition over $\mathcal B$ is
bounded by $\omega(\mathcal B)$. Actually, it can be easily seen that
if a graph has a tree-decomposition over a minor-closed class
$\mathcal B$ then it has a tree-decomposition over $\mathcal B$ of
adhesion at most $\omega(\mathcal B)-1$.

\sus{Path decompositions}
A \emph{path-decomposition} of a graph $G$ is a tree decomposition
where the underlying tree is a path. Of course we can always assume
that the path $P$ of a path decomposition $(P,(B_p)_{p\in P})$ has
vertex set $V^P=\{1,\ldots, m\}$, for some $m\in\mathbb N$, and that the
vertices occur on $P$ in their natural order (that is, we have
$i(i+1)\in E^P$ for $1\le i<m$).

\begin{lem}\label{attachlem}
  Let $G,H$ be graphs and $(\{1,\ldots,m\},(B_i)_{1\le i\le m})$ a
  path-decomposition of $H$ of width $k$. Let $x_1\ldots x_m$ be a
  path in $G$ such that $x_i\in B_{i}$ for $1\le i\le m$ and
  $V^{G}\cap V^{H}=\{x_1,\ldots,x_m\}$.  Then $ \tw(G\cup
  H)\le(\tw(G)+1)(k+1)-1.  $
\end{lem}

\proof
Let $(T,(C_t)_{t\in V^T})$ be a tree-decomposition of $G$. Then
$(T,(C'_t)_{t\in V^T})$ with  
\[
C'_t=C_t\cup\bigcup_{\substack{1\le i\le m,\\ x_i\in C_t}}B_{i}
\]
is a tree-decomposition of $G\cup H$.
\proofend

\section{Local tree-width}\label{sec:ltw}

The distance $d^G(x,y)$ between two vertices $x,y$ of a graph $G$ is
the length of the shortest path in $G$ from $x$ to $y$. For $r
\ge 1$ and $x \in G$ we define the \emph{$r$-neighborhood} around $x$ to be
$N^G_r(x) := \{ y \in V^G\mid d^G(x,y) \le r \}$. 

\begin{defn}
\begin{enumerate}
\item
The \emph{local tree-width} of a graph $G$ is the function $\ltw^G:\mathbb N\rightarrow\mathbb N$ defined by
\[
\ltw^G(r):=\max\big\{\tw(\langle N^G_r(x)\rangle)\bigmid x\in
V^G\big\}.
\]

\item
A class $\mathcal C$ of graphs has \emph{bounded local tree-width} if
  there is a function $f:\mathbb N\rightarrow\mathbb N$ such that
  $\ltw^{G}(r)\le f(r)$ for all $G\in\mathcal C$, $r\in\mathbb N$.

$\mathcal C$ has \emph{linear local tree-width} if
there is a $\lambda\in\mathbb R$ such that $\ltw^G(r)\le \lambda r$ for all
$G\in\mathcal C$, $r\in\mathbb N$.
\end{enumerate}
\end{defn}

\begin{ex}
  Let $G$ be a graph of tree-width at most $k$. Then $\ltw^G(r)\le k$ for all 
  $r\in\mathbb N$.
\end{ex}

\begin{ex}
  Let $G$ be a graph of valence at most $l$, for an $l\ge 1$. Then
  $\ltw^G(r)\le l(l-1)^{r-1}$ for all $r\in\mathbb N$.
\end{ex}

The planar graph algorithms due to Baker and others that we mentioned
in the introduction are based on the following result:

\begin{nsatz}{Proposition}{Bodlaender \cite{bod88}}
The class of planar graphs has linear local tree-width. More
precisely, for every planar graph $G$ and $r\ge 1$ we have
$\ltw^G(r)\le 3r$.
\end{nsatz}

In this paper, a \emph{surface} is a compact connected 2-manifold with
(possibly empty) boundary. The (orientable or non-orientable)
\emph{genus} of a surface $S$ is denoted by $g(S)$.
An \emph{embedding} of a graph $G$ in a
surface $S$ is a
mapping $\Pi$ that associates distinct points of $S$ with the vertices
of $G$ and internally disjoint simple curves in $S$ with the edges
of $G$ in such a way that a vertex $v$ is incident with an edge $e$ if, and
only if, $\Pi(v)$ is an endpoint of $\Pi(e)$.

\begin{nsatz}{Proposition}{Eppstein \cite{epp99}}\label{ltw-emb}
Let $S$ be a surface. Then the class of all graphs embeddable in $S$
has linear local tree-width. More precisely, there is a constant $c$
such that for all graphs $G$ embeddable in $S$ and for all $r\ge 0$ we
have $\ltw^G(r)\le c\cdot g(S)\cdot r$. 
\end{nsatz}

In the next subsection, we prove an extension of Proposition
\ref{ltw-emb} that forms the bases of our decomposition theorem for
graphs with excluded minors.

But before we do so, let me state another result due to Eppstein that
characterizes the minor closed classes of graphs of bounded
local tree-width. An \emph{apex graph} is a graph $G$ that has a
vertex $v\in V^G$ such that $G\setminus\{v\}$ is planar.

\begin{nsatz}{Theorem}{Eppstein \cite{epp95,epp99}}
Let $\mathcal C$ be a minor-closed class of graphs. Then $\mathcal C$
has bounded local tree-width if, and only if, $\mathcal C$ does not contain all
apex graphs.
\end{nsatz}

It is an interesting open problem whether there is a minor closed
class of graphs of bounded local tree-width that does not have linear
(or polynomially bounded) local tree-width.

\sus{Almost embeddable graphs}
Let $S$ be a surface with non-empty boundary. The boundary of $S$
consists of finitely many connected components
$C_1,\ldots,C_\kappa$, each of which is homeomorphic to the cycle $S^1$. 

We now define a graph $G$ to be \emph{almost embeddable} in $S$.
Roughly, this means that we can obtain $G$ from a graph $G_0$ embedded
in $S$ by attaching at most $\kappa$ graphs of path-width at most $\kappa$ to
$G_0$ along the boundary cycles $C_1,\ldots,C_\kappa$ in an orderly way.

This notion plays an important role in the structure theory of graphs
with excluded minors, to be outlined in the next subsection. 

\begin{defn}\label{almostembeddable}
Let $S$ be a surface with boundary cycles $C_1,\ldots,C_\kappa$. A graph $G$
is \emph{almost embeddable} in $S$ if there are
(possibly empty) subgraphs $G_0,\ldots,G_\kappa$ of $G$ such that
\begin{itemize} 
\item
$G=G_0\cup\ldots\cup G_\kappa$,
\item
$G_0$ has an embedding $\Pi$ in $S$, 
\item
$G_1,\ldots,G_\kappa$ are pairwise disjoint,
\item
for $1\le i\le \kappa$, $G_i$ has a path decomposition
$(\{1,\ldots,m_i\},(B^i_j)_{1\le j\le m_i})$ of width at most $\kappa$,
\item for $1\le i\le \kappa$ there are vertices
  $x^i_1,\ldots,x^i_{m_i}\in V^{G_0}$ such that $x^i_j\in B^i_j$
  for $1\le j\le m_i$ and $V^{G_0}\cap
  V^{G_i}=\{x^i_1,\ldots,x^i_{m_i}\}$,
\item for $1\le i\le \kappa$, we have $\Pi(V^{G_0})\cap
  C_i=\{\Pi(x^i_1),\ldots,\Pi(x^i_{m_i})\}$, and the points
  $\Pi(x^i_1),\ldots,\Pi(x^i_{m_i})$ appear on $C_i$ in this order
  (either if we walk clockwise or anti-clockwise).
\end{itemize}
\end{defn}

\begin{prop}\label{ltwalmostemb}
Let $S$ be a surface. Then the class of all graphs almost embeddable in $S$
has linear local tree-width.
\end{prop}

\begin{figure}[ht]
\begin{center}
\includegraphics*[height=3cm]{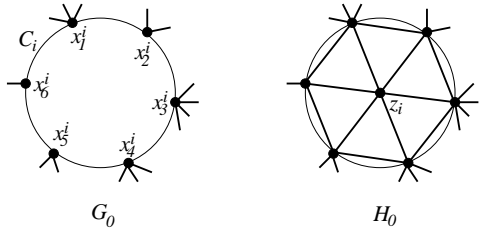}
\end{center}
\caption{From $G_0$ to $H_0$}\label{zifig}
\end{figure}

\proof Let $G$ be a graph that is almost embeddable in $S$.  We use
the notation of Definition \ref{almostembeddable}. Let $H_0$ be the
graph obtained from $G_0$ by adding new vertices
$z_1,\ldots,z_\kappa$, and edges $(z_i,x^i_j)$, $(x^i_j,x^i_{j+1})$,
and $(x^i_\kappa,x^i_1)$, for $1\le i\le \kappa$, $1\le j\le m_i$
(see Figure \ref{zifig}).  Clearly, $H_0$ is still embeddable in $S$.
For $1\le i\le \kappa$ we let $H_i:=H_0\cup G_1\cup\ldots\cup G_i$.

Let $\lambda\in\mathbb N$ such that for every graph $G$ embedabble in
$S$ and every $r\in\mathbb N$ we have $\ltw^G(r)\le\lambda r$ (such a
$\lambda$ exists by Theorem \ref{ltw-emb}). For $r\in\mathbb N$ we let
$f_0(r):=\lambda r$ and, for $i\in\mathbb
N$, we let $f_{i}(r):=(f_{i-1}(r+1)+1)(\kappa+1)-1$. Then $f_i$ is a
linear function for every $i\in\mathbb N$.

By induction on $i\ge 0$ we shall prove that for every $r\in\mathbb
N$ and $x\in V^{H_i}$ we have
\begin{equation}\label{atwind}
\tw\big(\langle N_r^{H_i}(x)\rangle\big)\le f_i(r).
\end{equation}
For $i=0$, this is immediate. So we assume that $i\ge 1$ and that we
have proved \eqref{atwind} for $i-1$.

For all $x\in H_i$, we either have $N_r^{H_i}(x)\subseteq H_{i-1}$, or
$N_r^{H_i}(x)\subseteq G_i$, or
$N_r^{H_i}\cap\{x^i_1,\ldots,x^i_{m_i}\}\neq\emptyset$.

If $N_r^{H_i}(x)\subseteq V^{H_{i-1}}$ then $\tw\big(\langle
N_{r}^{H_i}(x)\rangle^{H_i}\big)\le f_{i-1}(r)\le f_i(r)$.

If $x\in V^{H_{i-1}}$ and $N_r^{H_i}(x)\not\subseteq V^{ H_{i-1}}$, then
$N_{r-1}^{H_i}(x)\cap\{x^i_1,\ldots,x^i_{m_i}\}\neq\emptyset$. By the
construction of $H_0$, this implies $z_i\in N_r^{H_{i-1}}(x)$ and thus
$\{x^i_1,\ldots,x^i_{m_i}\}\subseteq N_{r+1}^{H_{i-1}}(x)$.

By Lemma \ref{attachlem} and the induction hypothesis we get
\begin{align*}
\tw\big(\langle N_{r}^{H_i}(x)\rangle^{H_i}\big)&\le
\tw\big(\langle N_{r+1}^{H_{i-1}}(x)\cup V^{G_i}\rangle^{H_i}\big)\\
&\le
( f_{i-1}(r+1)+1)(\kappa+1)-1= f_i(r).
\end{align*}

If $x\in V^{G_i}$, then 
$N_r^{H_i}(x)\cap V^{H_{i-1}}\subseteq N_{r+1}^{H_{i-1}}(z_i)$.
Thus by Lemma \ref{attachlem} and the induction hypothesis we have
\begin{align*}
\tw\big(\langle N_{r}^{H_i}(x)\rangle^{H_i}\big)&\le
\tw\big(\langle N_{r+1}^{H_{i-1}}(z_i)\cup V^{G_i}\rangle^{H_i}\big)\\
&\le
( f_{i-1}(r+1)+1)(\kappa+1)-1= f_i(r).
\end{align*}
\proofend

Note that the local tree-width of a graph is not minor-monotone (that
is, $H\preceq G$ does not imply $\ltw^H(r)\le\ltw^G(r)$ for all $r$). However, we do
have 
\begin{equation}\label{eq:ltw}
H\subseteq G\implies\ltw^H\le\ltw^G.
\end{equation}

\begin{prop}\label{ltwminalmostemb}
  Let $S$ be a surface. Then the class of all minors of graphs almost
  embeddable in $S$ has linear local tree-width.
\end{prop}

\proof Recall the proof of Proposition \ref{ltwalmostemb}. We use the
same notation here. Suppose $G'$ is a minor of $G$. We can assume that
$G'$ is a subgraph of a graph $G''$ obtained from $G$ only by
contracting edges. Because of \eqref{eq:ltw} we can even assume that $G'=G''$.

Let $X=\{x^i_j\mid 1\le i\le \kappa,1\le j\le m_i\}$. Contracting
edges with at least one endpoint not in $X$ is unproblematic, because
the resulting graph is still almost embeddable in $S$.

So we can further assume that $G'$ is obtained from $G$ by contracting
edges $e_1$, $\ldots$, $e_n$ with both endpoints in $X$. Let $H:=H_\kappa$ (the graph obtained
from $G$ by adding the vertices $z_i$ and corresponding edges as in Figure
\ref{zifig}). Let $H'$ be the graph obtained from $H$ by contracting
$e_1,\ldots,e_n$, and let $h:H'\preceq H$ witness these edge contractions.

\medskip
The key observation is that for all $x,y\in V^{H'}$ and $u\in h(x), v\in
h(y)$ we have
\begin{equation}\label{eq:twmin}
d^H(u,v)\le d^{H'}(x,y)+3\kappa-1
\end{equation}
(no matter how large $n$ is). To see this, let $P'$ be a shortest path
from $x$ to $y$ in $H'$. Let $P$ be a path from $u$ to $v$ in $H$ such
that $P'$ is obtained from $P$ by contracting the edges
$e_1\ldots,e_n$. Let us call such an edge an $(i,j)$-edge if it
connects a vertex in $\{x^i_1,\ldots,x^i_{m_i}\}$ with a vertex in
$\{x^{j}_1,\ldots,x^{j}_{m_j}\}$. Suppose that $P=w_1\ldots w_r$.  For
$1\le i\le \kappa$, let $w_s$ and $w_{t}$, where $1\le s\le t\le r$,
be the first and last vertex from $\{x^i_1,\ldots,x^i_{m_i}\}$ on $P$.
If $s<t$ we replace the interval $w_s\ldots w_t$ in $P$ by
$w_sz_iw_t$.  Doing this for $1\le i\le \kappa$ we obtain a new path
$Q$ from $u$ to $v$ in $H$. This path $Q$ contains no at most
$2\kappa$ edges that are not on $P$ and no $(i,i)$-edges. Furthermore,
for $1\le i<j\le n$ the number of $(i,j)$-edges on $Q$ is at most
$(\kappa-1)$. Because assume that $Q$ contains at least $\kappa$ such
edges. Then there would be a ``cycle'' $i=i_1,i_2,\ldots,i_l=i$ such
that for $1\le j<l$, $Q$ contains an $(i_j,i_{j+1})$-edge.  However,
this cycle would have been removed while transforming $P$ to $Q$.

Hence $\text{length}(Q)\le\text{length}(P')+3\kappa-1$, which proves
\eqref{eq:twmin}.

\medskip
\eqref{eq:twmin} implies that for all $r\ge 0$, $x\in V^{H'}$, and $u\in h(x)$ we have
\begin{equation}\label{eq:twmin2} 
\langle N_r^{H'}(x)\rangle\preceq\langle N^H_{r+3\kappa-1}(u)\rangle.
\end{equation}
To see this, let $y\in N^{H'}_r(x)$. Then for all $v\in h(y)$, by
\eqref{eq:twmin} we have $v\in N^H_{r+3\kappa-1}(u)$. Thus
$h(N^{H'}_r(x))\subseteq\Pow(N^{H}_{r+3\kappa-1}(u))$. Therefore the
restriction of $h$ to $N_r^{H'}(x)$ witnesses $\langle
N_r^{H'}(x)\rangle\preceq\langle N^H_{r+3\kappa-1}(u)\rangle$.
This proves \eqref{eq:twmin2}.

\medskip
By \eqref{atwind} and \eqref{eq:twmin2} we get $\tw(\langle N^{H'}_r(x)\rangle)\le
f_\kappa(r+3\kappa-1)$. The statement of the lemma follows.
\proofend

\section{Graphs with excluded minors}\label{sec:minclosed}
The following deep structure theorem for $K_n$-free graphs plays a
central role in the proof of the Graph Minor Theorem. For a surface
$S$ and $\mu\in\mathbb N$ we let $\mathcal A(S,\mu)$ be the class of
all graphs $G$ such that there is an $X\subseteq V^G$ with $\|X\|\le \mu$
such that $G\setminus X$ is almost embeddable in $S$.

\begin{nsatz}{Theorem}{Robertson and Seymour \cite{GMXVI}}\label{gm16}
  For every $n\in\mathbb N$ there exist $\mu\in\mathbb N$ and
  surfaces $S, S'$ such that all $K_n$-free graphs have a
  tree-decomposition over  $\mathcal A(S,\mu)\cup\mathcal A(S',\mu)$.
\end{nsatz}

Further details concerning this theorem can be found in 
\cite{dietho99,rtho99,GMXVI}.

For $\lambda,\mu\ge 0$ we
let
\begin{align*}
\mathcal L(\lambda)&:=\big\{ G\bigmid\forall H\preceq G\;\forall
r\ge0:\;\ltw^H(r)\le\lambda\cdot r\big\},\\
\mathcal L(\lambda,\mu)&:=\Big\{ G\Bigmid\exists X\subseteq
V^G:\big(\|X\|\le \mu\wedge G\setminus X\in\mathcal L(\lambda)\big)\Big\}.
\end{align*}
Note that $\mathcal L(\lambda,\mu)$ is minor closed and that $\omega(\mathcal L(\lambda,\mu))=\lambda+\mu+1$. Thus a
tree-decomposition over $\mathcal L(\lambda,\mu)$ has adhesion at most
$\lambda+\mu+1$. 

\begin{theo}\label{theo:minorclosed}
  Let $\mathcal C$ be a class of graphs with an excluded minor. Then
  there exist $\lambda,\mu\in\mathbb N$ such that all $G\in\mathcal C$ have a
  tree-decomposition over $\mathcal L(\lambda,\mu)$.
\end{theo}

\proof
This follows immediately from Theorem \ref{gm16} and Proposition \ref{ltwminalmostemb}.
\proofend

For algorithmic applications we have in mind, Theorem
\ref{theo:minorclosed} alone is not enough; we also have to compute
a tree-decomposition of a given graph over $\mathcal L(\lambda,\mu)$.
Fortunately, Robertson and Seymour have proved another deep result
that helps us with this task:

\begin{nsatz}{Theorem}{Robertson and Seymour \cite{GMXIII}}\label{gm13}
\sloppy
Every minor closed class of \mbox{graphs} has a
  polynomial time membership test.
\end{nsatz}

\begin{lem}\label{treedeclem}
Let $\mathcal C$ be a minor closed class of
graphs.

Then there is a polynomial time algorithm that computes, given
a graph $G$, a tree-decomposition of $G$ over $\mathcal C$, or rejects $G$ if no such tree-decomposition exists.
\end{lem}

\proof 
Note that the class $\mathcal T$ of all graphs that have a
tree-decomposition over $\mathcal C$ is minor closed. Thus by Theorem \ref{gm13} we have polynomial
time membership tests for both $\mathcal C$ and $\mathcal T$.

Without loss of generality, we can assume that $\mathcal C$ is not the
class of all graphs. Thus the clique number $\omega:=\omega(\mathcal C)$ is
finite. Recall that every
tree-decomposition over $\mathcal C$ has adhesion at most $\omega$.
Our algorithm uses the following observation to recursively construct a
tree-decomposition of the input graph $G$:
\begin{quote}\itshape
  $G\in\mathcal T$ if, and only if, $G\in\mathcal C$ or there is a set
  $X\subseteq V^G$ such that $|X|\le \omega$, $G\setminus  X$ has at least
  two connected components, and for all components $C$ of $G\setminus X$
  we have $\langle X\cup C\rangle^G\cup K_X\in\mathcal T$.
\end{quote}
We omit the details.
\proofend

In particular, we are going to apply this result to the minor closed
classes $\mathcal L(\lambda,\mu)$.

\section{Approximation algorithms}\label{sec:aa}

\sus{Optimization problems}\label{optprob}
An \emph{NP-optimization problem} is a tuple $(I,S,C,\opt)$,
consisting of a polynomial time decidable set $I$ of \emph{instances},
a mapping $S$ that associates a non-empty set $S(x)$ of
\emph{solutions} with each $x\in I$ such that the binary relation
$\{(x,y)\mid y\in S(x)\}$ is polynomial time computable and there is a
$k\in\mathbb N$ such that for all $x\in I$, $y\in S(x)$ we have
$\|y\|\le\|x\|^k$, a polynomial time computable \emph{cost} (or
\emph{value}) function $C:\{(x,y)\mid x\in I, y\in
S(x)\}\rightarrow\mathbb N$, and a \emph{goal} $\opt\in\{\min,\max\}$.

Given an $x\in I$,
we want to find a $y\in S(x)$ such that
\[
C(x,y)=\opt(x):=\opt\{C(x,z)\mid z\in S(x)\}.
\]

Let $x\in I$ and $\epsilon>0$. A solution $y\in S(x)$ for $x$ is
\emph{$\epsilon$-close} if
\[
(1-\epsilon)\opt(x)\le C(x,y)\le (1+\epsilon)\opt(x).
\]
A \emph{polynomial time approximation scheme (PTAS)} for
$(I,S,C,\opt)$ is a uniform family $(A_\epsilon)_{\epsilon>0}$ of
approximation algorithms, where $A_\epsilon$ is a polynomial time
algorithm that, given an $x\in I$,  computes
an $\epsilon$-close solution for $x$ in polynomial time. Uniformity
means that there is an algorithm that, given $\epsilon$, computes
$A_\epsilon$.

\sus{The levels of graphs of bounded local tree-width}
For  graph $G$, a vertex $v\in V^G$, and integers $j\ge i\ge 0$ we let
\[
L^G_v[i,j]:=\{w\in V^G\mid i\le d^G(v,w)\le j\}.
\]
To keep the notation uniform, we are actually going to write
$L_v^G[i,j]$ for arbitrary $i,j\in\mathbb Z$, with the understanding
that $L_v^G[i,j]:=\emptyset$ for $i>j$ and $L_v^G[i,j]:=L_v^G[0,j]$
for $i\le 0$.

\begin{lem}\label{lem:level}
  Let $\lambda\in\mathbb N$. Then for all $G\in\mathcal L(\lambda)$,
  $v\in V^G$, and $i,j\in\mathbb Z$ with $i\le j$ we have $\tw\big(\langle
  L_v^G[i,j]\rangle\big)\le\lambda\cdot(j-i+1)$.
\end{lem}

\proof
First note that
$L_v^G[1,j]\subseteq L_{v}^G[0,j]=N_j^G(v)$, thus the claim holds for $i\le
1$.  For $i\ge2$, consider the minor $H$ of $G$ obtained by contracting the
connected subgraph $\langle L_{v}^G[0,i-1]\rangle$ to a single vertex $v'$.
Then we have $L_{v}^G[i,j]\subseteq N_{j-i+1}^{H}(v')$, and the claim follows.
\proofend

\sus{Minimum vertex cover}
Instances of \minvc\ are graphs $G$, solutions are sets $X\subseteq V^G$ such
that for every edge $vw\in E^G$ either $v\in X$ or $w\in X$ (such sets $X$
are called \emph{vertex covers}), the cost function is defined by
$C(G,X):=|X|$, and the goal is $\min$.

\begin{nsatz}{Lemma}{\cite{arnlagses91}}\label{lem:als}
For every $k\ge 1$, the restriction of\/ \minvc\/ to instances of
tree-width at most $k$ is solvable in linear time.
\end{nsatz}

\begin{theo}\label{theo:minvc}
Let $\mathcal C$ be a class of graphs with an excluded minor. Then the
restriction of \minvc\ to instances in $\mathcal C$ has a PTAS.
\end{theo}

\proof Applying Theorem \ref{theo:minorclosed}, we choose $\lambda,\mu\in\mathbb N$ such that
every $G\in\mathcal C$ has a tree-decomposition over $\mathcal
L(\lambda,\mu)$. Let $\epsilon>0$; we shall describe a polynomial time
algorithm that, given a graph $G\in\mathcal C$, computes an
$\epsilon$-close solution for \minvc\ on $G$. Uniformity will be
clear from our description. Let $k=\lceil\frac{1}{\epsilon}\rceil$ and
note that $\frac{k+1}{k}\le(1+\epsilon)$.

\medskip
In a first step, let us prove that the restriction of \minvc\ to instances in
$\mathcal L(\lambda)$ has a PTAS.

Let $G\in\mathcal L(\lambda)$ and $v\in V^G$ arbitrary.
For $1\le i\le k$ and $j\ge 0$ we let
$L_{ij}:=L_v^G[(j-1)k+i,jk+i]$. By Lemma \ref{lem:level}, $\tw(\langle
L_{ij}\rangle)\le\lambda(k+1)$. 

For $1\le i\le k$, $j\ge 0$
let $X_{ij}$ be a minimal vertex cover of $\langle L_{ij}\rangle$. We let $X_i:=\bigcup_{j\ge0}X_{ij}$. Then $X_i$ is a
vertex cover of $G$. Let $X_{\min}$ be a minimal vertex cover for
$G$. We have 
$|X_{ij}|\le |X_{\min}\cap L_{ij}|$, because $X_{\min}\cap L_{ij}$ is also a vertex cover of
$\langle L_{ij}\rangle$. Hence
\[
\sum_{i=1}^k|X_i|\le\sum_{i=1}^k\sum_{j\ge0}|X_{ij}|
\le\sum_{i=1}^k\sum_{j\ge0}|L_{ij}\cap X_{\min}|\le(k+1)|X_{\min}|.
\]
The last inequality follows from the fact that every $v\in V^G$ is
contained in at most $(k+1)$ (successive) sets $L_{ij}$.

Choose $m,1\le m\le k$ such that $|X_{m}|=\min\{|X_1|,\ldots,|X_k|\}$. Then
\[
|X_{m}|\le\frac{k+1}{k}|X_{\min}|\le(1+\epsilon)|X_{\min}|. 
\]
Since the $X_{ij}$ can be computed in polynomial time by Lemma
\ref{lem:als}, $X_{m}$ can also be computed in polynomial time. 

\medskip
In a second step, we show how to extend this approximation algorithm to classes $\mathcal
L(\lambda,\mu)$ for $\lambda,\mu\ge 0$. Let $G\in\mathcal
L(\lambda,\mu)$ and $U\subseteq V^G$ such that $|U|\le \mu$ and
$H:=G\setminus U\in\mathcal L(\lambda,0)$. The following extension of Lemma \ref{lem:als} can be proved
by standard dynamic programming techniques (cf.\ \cite{arnlagses91}):

\begin{lem}\label{lem:als2}
  For every $k\ge 0$, the following problem can be solved in linear
  time: Given a graph $G$, a subset $U\subseteq V^G$ such that
  $\tw(G\setminus U)\le k$, and a subset $Y\subseteq U$, compute a set
  $X\subseteq V^G\setminus U$ of minimal order such that $X\cup Y$ is a
  vertex cover of $G$, if such a set exists, or reject otherwise.
\end{lem}

For every $Y\subseteq U$ we shall compute an
$X(Y)\in\text{Pow}(V^G\setminus U)\cup\{\bot\}$ such that either $X(Y)\cup Y$ is a vertex cover of
$G$ and 
\[
|X(Y)|\le(1+\epsilon)\min\{|X|\mid X\subseteq V^G\setminus U, X\cup Y\text{ vertex cover of $G$}\},
\]
or $X(Y):=\bot$ if no such $X(Y)$ exists. Using Lemma \ref{lem:als2}
instead of Lemma \ref{lem:als}, we can do this analogously to the
first step.

Then we choose a $Y_0\subseteq U$ such that $|X(Y_0)\cup Y_0|$ is
minimal. Here we define $\bot\cup Z:=\bot$ for all $Z$ and
$|\bot|:=\infty$. Then clearly $X(Y_0)\cup Y_0$ is an $\epsilon$-close
solution for \minvc\ on $G$. Moreover, since $|U|\le\mu$, there are
at most $2^\mu$ sets $Y\subseteq U$, so $X(Y_0)\cup Y_0$ can be
computed in polynomial time (remember that $\mu$ is a constant only
depending on the class $\mathcal C$).

\medskip
In the third step, we extend our PTAS to graphs that have a
tree-decomposition over $\mathcal L(\lambda,\mu)$, i.e.\ to all graphs
in $\mathcal C$.

So let $G$ be such a graph. We first compute a tree-decomposition
$(T,(B_t)_{t\in V^T})$ of $G$ over $\mathcal L(\lambda,\mu)$. Remember
that by Lemma \ref{treedeclem}, this is possible in polynomial time.
Recall that $r^T$ denotes the root of $T$ and that, for every $t\in V^T$
with parent $u$, we let $A_t=B_t\cap B_u$. For every $t\in V^T$, we let
$S_t$ be the subtree of $T$ with root $t$, that is, the subtree with
vertex set $\{s\mid t\text{ occurs on the path from $s$ to $r^T$}\}$. We let
$C_t:=\bigcup_{s\in S_t}B_t$.

Inductively from the leaves to the root, for every
node $t\in V^T$ and for every $Y\subseteq A_t$
we compute an $X(t,Y)\in\text{Pow}(C_t\setminus
A_t)\cup\{\bot\}$ such that either
$X(t,Y)\cup Y$
is a vertex cover of $\langle C_t\rangle$ and 
\[
|X(t,Y)|\le(1+\epsilon)\min\{|X|\mid X\cup Y\text{ vertex cover of $\langle
  C_t\rangle$}\},
\]
or $X(t,Y):=\bot$ if no such vertex set exists.
Since a tree-decomposition over $\mathcal L(\lambda,\mu)$ has adhesion
at most $\lambda+\mu+1$ we have $|A_t|\le\lambda+\mu+1$, thus for every
$t\in V^T$ we have to compute at most $2^{\lambda+\mu+1}$ sets
$X(t,Y)$. For the root $r^T$ we have $A_{r^T}=\emptyset$, so
$X(r^T,\emptyset)$ is an $\epsilon$-close solution for \minvc\ on $G$.

Suppose that $t\in V^T$ and that for every child $t'$ of $T$ we have
already computed the family $X(t',\cdot)$. Let $U\subseteq B_t$ such
that $|U|\le \mu$ and $[B_t]\setminus U\in\mathcal L(\lambda)$. Let
$W:=U\cup A_t$ and let $Z\subseteq W$.  Let
$X_{\min}(Z)\in\text{Pow}(C_t\setminus W)\cup\{\bot\}$ be a vertex set
of minimal order such that $X_{\min}(Z)\cup Z$ is a vertex cover of
$\langle C_t\rangle$, or $X(Z):=\bot$ if no such vertex set exists.

We show how to compute an $X(Z)\in\text{Pow}(C_t\setminus
W)\cup\{\bot\}$ such that $X(Z)\cup Z$ is a vertex cover of $\langle
C_t\rangle$ and $|X(Z)|\le(1+\epsilon)|X_{\min}(Z)|$, if
$X_{\min}(Z)\neq\bot$, or $X(Z)=\bot$ otherwise. Then for every
$Y\subseteq A_t$ we choose a $Z\subseteq W$ such that $Y\subseteq Z$
with minimal $|X(Z)\cup(Z\setminus Y)|$ (among all $Z\supseteq Y$) and
let $X(t,Y):=X(Z)$.  Note that, since $|U|\le\mu$, for every $Y$ we
have to compute at most $2^\mu$ sets $X(Z)$ to determine $X(t,Y)$.

So let us fix a $Z\subseteq W$; we show how to
compute $X(Z)$ in polynomial time.

If $W=B_t$ we let $X(Z):=\bigcup_{t'\text{ child of }t}X(t',A_{t'}\cap
Z)$.

Otherwise, we choose an arbitrary $v\in B_t\setminus W$. For $1\le
i\le k$ and $j\ge 0$ we let $L_{ij}:=L_v^{[B_t]\setminus
  W}[(j-1)k+i,jk+i]$. Then $\tw(\langle
L_{ij}\rangle)\le\lambda(k+1)$. For $1\le i\le k$ and every child $t'$
of $t$ there is at least one $j\ge 0$ such that $A_{t'}\setminus
W\subseteq L_{ij}$, because $A_{t'}$ induces a clique in $[B_t]$. Let
$j^*(i,t')$ be the least such $j$ and
$L_{ij}^*:=L_{ij}\cup\bigcup_{\substack{t'\text{ child of
      }t\\j^*(i,t')=j}}C_{t'}\setminus A_{t'}$.

For every $X\subseteq L_{ij}$ we let
\[
X^*:=X\cup\bigcup_{\substack{t'\text{ child of
      }t\\j^*(i,t')=j}}X(t',(X\cup Z)\cap A_{t'})
\]
We compute an $X_{ij}\subseteq L_{ij}$ with minimal $|X_{ij}^*|$ such
that $X_{ij}\cup Z$ is a vertex cover of $\langle L_{ij}\cup W\rangle$
if such a vertex cover exists, and $X_{ij}=\bot$ otherwise. The usual
dynamic programming techniques on graphs of bounded tree-width show
that each $X_{ij}$ can be computed in linear time if the numbers
$|X(t',Y)|$ for the children $t'$ of $t$ are given (cf.\ Lemmas
\ref{lem:als} and \ref{lem:als2} and \cite{arnlagses91}). It is
important here that every $A_{t'}\setminus W$ is a clique in $\langle
L_{ij}\rangle$ and thus by Lemma \ref{lem:treedec}(1) completely
contained in a block of every tree-decomposition of $\langle
L_{ij}\rangle$. 
 
We let $X_{i}:=\bigcup_{j\ge0}X_{ij}$ and $X_{i}^*:=\bigcup_{j\ge0}X_{ij}^*$.
Then $X_i^*\cup Z$ is a vertex cover of $\langle C_t\rangle$, if such a
vertex cover exists, and $X_{i}=\bot$ otherwise. We choose
an $i,1\le i\le k$, such that $|X_i^*|=\min\{|X_1^*|,\ldots,|X_k^*|\}$ and
let $X(Z):=X_i^*$. Then $X(Z)$ can be computed in polynomial time.

Recall that $X_{\min}:=X_{\min}(Z)\subseteq C_t\setminus W$ is a
vertex set of minimal order such that $X_{\min}\cup Z$ is a vertex
cover of $\langle C_t\rangle$, if such a vertex cover exists, and
$X_{\min}=\bot$ otherwise. It remains to prove that
$|X(Z)|\le(1+\epsilon)|X_{\min}|$.

Recall that for every child $t'$ of $t$ we have 
\[
|X(t',(X_{\min}\cup Z)\cap A_{t'})|\le(1+\epsilon)|X_{\min}\cap
C_{t'}\setminus A_{t'}|.
\]
Our construction of the $X_{ij}$ and $X_{ij}^*$ guarantees that for
$1\le i\le k, j\ge0$ we have 
\[
|X_{ij}^*|\le |X_{\min}\cap L_{ij}|+\sum_{\substack{t'\text{ child of
        }t\\j^*(i,t')=j}}|X(t',(X_{\min}\cup Z)\cap A_{t'})|.
\]

Then
\begin{align*}
  k|X(Z)|\le&\sum_{i=1}^k|X_i^*|\\
  =&\sum_{i=1}^k\sum_{j\ge0}|X_{ij}^*|\\
  \le&\sum_{i=1}^k\sum_{j\ge0}\Big(|X_{\min}\cap
  L_{ij}|+\sum_{\substack{t'\text{ child of
        }t\\j^*(i,t')=j}}|X(t',(X_{\min}\cup Z)\cap A_{t'})|\Big)\\
  \le&\sum_{i=1}^k\sum_{j\ge0}\Big(|X_{\min}\cap
  L_{ij}|+\sum_{\substack{t'\text{ child of
        }t\\j^*(i,t')=j}}(1+\epsilon)|X_{\min}\cap C_{t'}\setminus
  A_{t'}|\Big)\\
\le&(k+1)|X_{\min}\cap B_t|+k(1+\epsilon)|X_{\min}\cap C_t\setminus B_t|.
\end{align*}
This implies $|X(Z)|\le(1+\epsilon)X_{\min}$.
\proofend

\sus{Minimum dominating set}
Instances of \minds\ are graphs $G$, solutions are sets $X\subseteq
V^G$ such that for every $v\in V^G\setminus X$ there is a $w\in X$
such that $vw\in E^G$ (such sets $X$ are called \emph{dominating
  sets}), the cost function is defined by $C(G,X):=|X|$, and the goal
is $\min$.

\begin{theo}
Let $\mathcal C$ be a class of graphs with an excluded minor. Then the
restriction of \minds\ to instances in $\mathcal C$ has a PTAS.
\end{theo}

\proof
We proceed very similarly to the proof of Theorem \ref{theo:minvc},
the analogous result for  \minvc. Let $\lambda,\mu\in\mathbb N$
such that every graph in $\mathcal C$ has a tree-decomposition over
$\mathcal L(\lambda,\mu)$. Let $\epsilon>0$ and
$k:=\lceil\frac{2}{\epsilon}\rceil$.

Again, in the first step we consider the restriction of the problem to
input graphs from $\mathcal L(\lambda)$. Given such a graph $G$, we
choose an arbitrary $v\in V^G$. For $1\le i\le k$ and $j\ge 0$ we let
$L_{ij}:=L^G_v[(j-1)k+i-1,jk+i]$. Then $\tw(\langle
L_{ij}\rangle)\le\lambda(k+2)$. Note that $L_{ij}$ and $L_{i(j+1)}$
overlap in two consecutive rows, which is different from the proof of
Theorem \ref{theo:minvc}. The \emph{interior} of $L_{ij}$ is the set $L_{ij}^\circ:=L^G_v[(j-1)k+i,jk+i-1]$.

For $1\le i\le k,j\ge 0$ we let $X_{ij}\subseteq L_{ij}$ be a vertex
set of minimal order with the following property:
\begin{itemize}
\item[$(*)$]
For every $w\in L_{ij}^\circ\setminus X_{ij}$ there is a $x\in X_{ij}$
such that $(w,x)\in E^{G}$. 
\end{itemize}
Then for $1\le i\le k$ the set
$X_i:=\bigcup_{j\ge 0}X_{ij}$ is a dominating set of $G$. Let $m$ be such that $|X_{m}|=\min\{|X_1|,\ldots,|X_k|\}$. Computing 
$X_{m}$ amounts to solving a variant of \minds\ on instances of
tree-width at most $\lambda(k+2)$; using the usual dynamic programming
techniques, this can be done in linear time.

Since for every dominating set $X$ of $G$ the set $X\cap L_{ij}$ has
property $(*)$ we have $X_{ij}\le X\cap L_{ij}$. Using this, we can
argue as in the proof of Theorem \ref{theo:minvc} to show that
$X_{m}$ is an $\epsilon$-close solution.

Adapting the second and third step of the proof of Theorem
\ref{theo:minvc}, it is straightforward to extend this algorithm to
arbitrary input graphs in $\mathcal C$.
\proofend

\sus{Maximum independent set} Instances of \maxind\ are graphs $G$,
solutions are sets $X\subseteq V^G$ such that for all $v,w\in X$ we
have $vw\not\in E^G$ (such sets $X$ are called \emph{independent
  sets}), the cost function is defined by $C(G,X):=|X|$, and the goal
is $\max$.

\begin{theo}
Let $\mathcal C$ be a class of graphs with an excluded minor. Then the
restriction of \maxind\ to instances in $\mathcal C$ has a PTAS.
\end{theo}

\proof
Again we proceed similarly to the proof of Theorem \ref{theo:minvc}. 
Let $\lambda,\mu\in\mathbb N$
such that every graph in $\mathcal C$ has a tree-decomposition over
$\mathcal L(\lambda,\mu)$. Let $\epsilon>0$ and
$k=\lceil\frac{1}{\epsilon}\rceil$. 

We describe how to treat input graphs in $\mathcal L(\lambda)$. Following the lines of the proof of Theorem \ref{theo:minvc}, the
extension to arbitrary $G\in\mathcal C$ is straightforward. Let
$G\in\mathcal L(\lambda)$ and $v\in V^G$. For $1\le i\le k$ and
$j\ge 0$ we let $L_{ij}:=L^G_v[(j-1)k+i,jk+i-2]$. Then $\tw(\langle
L_{ij}\rangle)\le\lambda(k-1)$. Note that there are no edges between
$L_{ij}$ and $L_{i(j+1)}$.

For $1\le i\le k,j\ge 0$ we let $X_{ij}$ be a maximal independent set
of $\langle L_{ij}\rangle$. Then $X_i:=\bigcup_{j\ge 0} X_{ij}$ is an
independent set of $G$. Let $1\le m\le k$ such that
$|X_{m}|=\max\{|X_1|,\ldots,|X_k|\}$. Since the restriction of
\maxind\ to graphs of bounded tree-width is solvable in linear time,
such an $X_m$ can be computed in linear time.

Let $X_{\max}$ be a maximum independent set
of $G$. Then for $1\le i\le k$, $j\ge 0$  we have $|X_{ij}|\ge
|X_{\max}\cap L_{ij}|$.  Thus
\[
k|X_{m}|\ge\sum_{i=1}^k|X_i|=\sum_{i=1}^k\sum_{j\ge
  0}|X_{ij}|\ge\sum_{i=1}^k\sum_{j\ge 0}|X_{\max}\cap L_{ij}|\ge(k-1)|X_{\max}|,
\]
which implies that
$X_m\ge\frac{k-1}{k}|X_{\max}|\ge(1-\epsilon)|X_{\max}|$.
\proofend

\sus{Other problems}
Our approach can be used to find polynomial time approximation schemes
for the restrictions of a number of other problems to classes of
graphs with excluded minors, in particular for the other problems
considered by Baker \cite{bak94}.
I leave it to the reader to work out the details.

\section{Other applications of Theorem \ref{theo:minvc}}\label{sec:appl}
\sus{The tree-width of $K_n$-free graphs}\label{sec:ast}
We re-prove a theorem of Alon, Seymour, and Thomas \cite{aloseytho90a} that the
tree-width of a $K_n$-free graph $G$ is $O(\sqrt{|G|})$. This is joint
work with Reinhard Diestel and Daniela K\"uhn.

\begin{lem}
Let $\lambda\in\mathbb N$ and $G\in\mathcal L(\lambda)$. Then
$\tw(G)\le3\sqrt{\lambda |G|}$. 
\end{lem}

\proof
Let $v\in V^G$ arbitrary and, for $i\ge 0$, $L_i:=\{w\in V^G\mid d^G(v,w)=i\}$. 
Let $m$ be maximal such that $L_m$ is non-empty. We subdivide $\{1,\ldots,m\}$
into intervals $I_1,J_1,I_2,\ldots,J_{l-1},I_l,J_l$ such that for
$1\le i\le l$ we have
\begin{itemize}
\item
$|L_j|\le\sqrt{\lambda\cdot |G|}$ for all $j\in I_i$,
\item
$|L_j|>\sqrt{\lambda\cdot |G|}$ for all $j\in J_i$.
\end{itemize}
Then
$
\tw(\langle \bigcup_{j\in I_i}L_j\rangle)\le 2\sqrt{\lambda\cdot |G|}
$
and 
$
\tw(\langle \bigcup_{j\in J_i}L_j\rangle)\le \sqrt{\lambda\cdot |G|}
$
(because the length of $J_i$ is at most $\sqrt{\frac{|G|}{\lambda}}$).
We can glue the decompositions together by adding to every block of a
tree-decomposition of $J_i$ the last level of
the previous $I_{i}$ and the first level of the next $I_{i+1}$ and obtain $\tw(G)\le 3\sqrt{\lambda\cdot |G|}$.
\proofend

\begin{cor}
Let $\lambda,\mu\in\mathbb N$ and $G\in\mathcal L(\lambda,\mu)$. Then
$\tw(G)\le3\sqrt{\lambda |G|}+\mu$. 
\end{cor}

\begin{cor}
Let $G$ be $K_n$-free. Then $\tw(G)\le O(\sqrt{|G|})$.
\end{cor}

\sus{Deciding first-order properties} In \cite{flugro99} we give
another algorithmic application of Theorem \ref{theo:minorclosed}. We
show that for every class $\mathcal C$ of graphs with an excluded
minor there is a constant $c>0$ such that for every property of graphs
that is definable in first order logic there is an $O(|G|^c)$-algorithm
deciding whether a given graph $G\in\mathcal C$ has this property.

For example, this implies that for every class $\mathcal C$ with an
excluded minor there is a constant $c$ such that for every graph $H$
there is an $O(|G|^c)$-algorithm testing whether a given graph $G\in\mathcal C$
has a subgraph isomorphic to $H$.

\section{Further research}
We have never specified the exponents and coefficients of the
polynomials bounding the running times of our algorithms; they seem to
be enormous. So our algorithms are only of theoretical interest. The
first important step towards improving the algorithms would be a
practically applicable algorithm for computing tree-decompositions of
graphs of small tree-width.  On the graph theoretic side, it would
probably help to prove Theorem \ref{theo:minorclosed} directly without
using Robertson's and Seymour's Theorem \ref{gm16}.

\medskip
The traveling salesman problem is another optimization problem that
has a PTAS on planar graphs \cite{grikoupap95,arogrikarkleiwol98}. It
would be interesting to see if this problem has a PTAS on class of
graphs with an excluded minor. 

\subsection*{Acknowledgements}
I thank Reinhard Diestel and J\"org Flum for helful
  comments on earlier versions of this paper.

\end{document}